\documentstyle{amsppt}
\magnification=1200

\title On The Complemented Subspaces of $X_{p}$
\endtitle
                                     
\topmatter
\author Dale E. Alspach
\thanks Research supported in part by NSF grant 
DMS 890237.
\endthanks
\endauthor
\address 
Department of Mathematics
Oklahoma State University
Stillwater, OK 74078-0613
\endaddress
\abstract In this paper we prove some results related to the
problem of isomorphically classifying the complemented subspaces of
$X_{p}$. We characterize the complemented subspaces of $X_{p}$  which
are isomorphic to $X_{p}$  by showing that such a space must contain
a canonical complemented subspace isomorphic to $X_{p}.$ We also give
some characterizations of complemented subspaces of $X_{p}$ isomorphic
to $\ell_{p}\oplus \ell_{2}.$
\endabstract

\endtopmatter
\document

\heading 0. Introduction \endheading

     Rosenthal \cite{R} introduced the ${\Cal L }_{p}$  space $X_{p}$  in 1971.
Among its interesting
properties are that it contains and is contained in isomorphs of
$\ell_{p}\oplus \ell_{2},$ but
is not isomorphic to a complemented subspace of $\ell_{p}\oplus  \ell_{2}.$ 
These properties
have made $X_{p}$  rather resistant to standard approaches to classifying its
complemented subspaces. For example it was first proved that $X_{p}$  
was primary
in \cite{JO2} where the device of simultaneously $L_{p},\; L_{2}$  bounded 
operators was
employed to prove a version of the decomposition method for $X_{p}.$  The 
problem
really is that the $\ell_{p}$  and $\ell_{2}$  structures in $X_{p}$  are 
mixed in a much different
way than they are in $\ell_{p}\oplus \ell_{2}$  or $\left (\sum \ell_{2}
\right )_{p}.$  Let us also recall that $X_{p}$  or really
the technique for building $X_{p}$  is the central device used to construct an
uncountable number of separable ${\Cal L }_{p}$  spaces, \cite{BSR}. 
Thus a better understanding
of $X_{p}$  is critical for the study of the complemented subspaces of $L_{p}.$

     In order to state precisely our results we need to introduce some special
notation. Throughout this paper $w = (w_{n} )$ will be a sequence of 
positive real
numbers and $2 < p < \infty .$  As usual $X_{p,w}$    is the completion of
$\{ (a_{i} ):i \in \Bbb N ,\; a_{i}  \neq 0 \text{ for finitely many }
i\} $ with the 
norm $||(a_{i} )|| =
\max \{ |(a_{i} )|_{p},\; |(a_{i} )|_{2,w}\} $ where 
$ |(a_{i} )|_{p}  = [
\sum_{i} |a_{i} |^{p} ]^{1/p}$ and $|(a_{i} )|
_{2,w}=[
\sum_{i}|a_{i}|^{2}w^{2}_{i}]^{1/2}.$
    The Rosenthal space 
$X_{p}$  is $X_{p,w}$    where $w = (w_{i} )$ is such that for
every $\epsilon > 0,\; \sum\limits_{w_{i}<\epsilon }w^{2p/
(p-2)}_{i} = \infty .$  Throughout this paper we will always
consider $X_{p}$  to be the subspace of $(\ell_{p}\oplus \ell_{2} )_{\infty }$ 
 spanned by $(\delta_{n}  + w_{n} \gamma_{n} )$ where $(\delta_{n} ) $
and $(\gamma_{n} )$ are the usual unit vector bases of $\ell_{p}$  and 
$\ell_{2} ,$ respectively.  If $E\subset \Bbb N    $, the
symbol $\omega (E)$ will be used to denote  $\sum_{n\in 
E}w^{2p/(p-2)}_{n}$         which occurs frequently in
computations in $X_{p} .$  We will also need the ratio of 2-norm and
$p$-norm, 
$ |x|_{2} /|x|_{p} $,
which we will denote by $r(x)$. 
 For a subspace $Y$ of 
$X_{p} ,$ define $r(Y) =\sup \{ r(y): y \in  Y\} $ and $h(Y) = \inf \{ r(y):y
\in Y\} .$   In defining functionals on $X_{p}$  it
is convenient to use the inner product induced by the norm $|\cdot|_{2}.$  Thus
$<(x_{n} ),\; (y_{n} )> = \sum  x_{n} y_{n} w^{2}_{n}  .$

     We will use standard Banach space notation and terminology as may be
found in \cite{LT}. Here subspace will mean infinite dimensional closed subspace
unless otherwise noted.  The properties of $X_{p}$  can be found in 
\cite{LT,4d} or in
the original paper of Rosenthal \cite{R}.

\newpage
\heading 1. Complemented subspaces of $X_{p}$  which contain $X_{p}$  
complemented \endheading

     In this section we will show that a complemented subspace of $X_{p}$  which
contains a complemented subspace isomorphic to $X_{p}$  contains a canonical
complemented copy of $X_{p}.$  In \cite{R} Rosenthal showed that 
there were nice block
bases of the usual basis of $X_{p}$  with complemented closed linear span. 
The point
of his construction was to make sure that the coordinate functionals of the
block basis could be chosen to be bounded in both the $p$ and 2 norms.
Explicitly, if for each $j \in \Bbb N$    ,$$ 
     y_{j}=\sum^{k_{j+1}}_{n=k_{j}+1}w_{n}^{2/(p-2)}e_{n}$$ 
where $(e_{n} )$ is the natural basis for $X_{p}$  and $(e^{*}_{n})$ is the 
corresponding sequence of biorthogonal functionals, then $$ 
         |y_{j}|_{2}=\biggl [\sum^{k_{j+1}}_{n=k_{j}+1}w^{2p/(p-2)}_{n}
\biggr ]^{1/2} \text{ and }|y_{j}|_{p}=\biggl [
\sum^{k_{j+1}}_{n=k_{j}+1}w_{n}^{2p/(p-2)}\biggr ]^{1/p}.$$
Let $$
y^{*}_{j}(x)=|y_{j}|^{-2}_{2}\sum^{k_{j+1}}_{n=k_{j}+1}w_{n}
^{2/(p-2)}w^{2}_{n}e^{*}_{n}(x)=|y_{j}|^{-2}_{2}<y_{j},x>.$$
By applying H\"{o}lder's inequality first in $\ell_{2}$  and then in
$\ell_{p/2}$    we see that
$$\align
|y^{*}_{j}(x)| &= |y_{j} |^{-2}_{2} \biggl |
\sum^{k_{j+1}}_{n=k_{j}+1}w^{2/(p-2)}_{n}w^{2}_{n}e^{*}_{n}(x)\biggr|
\\
&\leq |y_{j}|^{-2}_{2}\biggl[
\sum^{k_{j+1}}_{n=k_{j}+1}w_{n}^{4/(p-2)}w^{2}_{n}\biggr]^{1/2}\biggl[
\sum^{k_{j+1}}_{n=k_{j}+1}w^{2}_{n}|e^{*}_{n}(x)|^{2}\biggr]^{1/2}\\
\split
&\leq |y_{j}|^{-2}_{2}\biggl[
\sum^{k_{j+1}}_{n=k_{j}+1}w_{n}^{2p/(p-2)}\biggr]^{1/2}\min \left\{
\biggl[\sum^{k_{j+1}}_{n=k_{j}+1}w^{2}_{n}|e^{*}_{n}(x)|^{2}\biggr]
^{1/2},\right.\\ 
&\qquad\qquad\qquad \left.\biggl[\sum^{k_{j+1}}_{n=k_{j}+1}w_{n}^{2p/(p-2)}
\biggr]^{(p-2)/2p}
\biggl[
\sum^{k_{j+1}}_{n=k_{j}+1}|e^{*}_{n}(x)|^{p}\biggr]
^{1/p}\right\} \endsplit\\
&=\min \left\{ |x_{\left| [k_{j}+1,
k_{j+1}]\right.} |_{2}|y_{j}|^{-1}_{2},\; |x_{\left| [k_{j}+1,k_{j+1}]\right.  }
|_{p}|y_{j}|_{p}^{-1}\right\} .
\endalign$$

Thus $$
|y_{j}^{*}(x)y_{j}|_{2}\leq |x_{| [k_{j}+1,k_{j+1}]}|_{2}
\text{ and }|y^{*}_{j}(x) y_{j}|_{p}\leq |x_{\left|
[k_{j}+1,k_{j+1}\right. ]}|_{p}.$$
The important point is that this computation works because
$$
\align
|y_{j}|_{2} &=\biggl[
\sum_{n=k_{j}+1}^{k_{j+1}}w_{n}^{2p/(p-2)}\biggr]^{(p-2)/2p}
\biggl[\sum^{k_{j+1}}_{n=k_{j}+1}w_{n}^{2p/(p-2)}\biggr]^{1/p}\\
&=\biggl[\sum_{n=k_{j}+1}^{k_{j+1}}w_{n}^{2p/(p-2)}\biggr]
^{(p-2)/2p}|y_{j}|_{p}.
\endalign$$

{\smc Remark 1.1:} The choice of coefficients for $y_{j}$  has the following 
geometric
motivation. These coefficients give the maximum for the ratio $|\cdot |_{2}  / 
|\cdot |_{p}$  for
elements supported on $[k_{j} +1,\; k_{j+1}].$  Thus if $x$ is supported on 
$[k_{j} +1,\;  k_{j+1}]$ and $x \in B_{\ell_{p}}$   (the unit ball of
$\ell_{p}$ ) then $(|y_{j}|_{p}/|y_{j}|_{2})x\in B_{|\cdot |_{2}}$ and there is
exactly one $x$ so that the multiple has 2-norm one.
\bigpagebreak

     If we replace $y_{j}$  by any $z_{j}$  with the same support which 
satisfies
$$
c|z_{j}|_{2}\geq [
\sum^{k_{j+1}}_{n=k_{j}+1}w_{n}^{2p/(p-2)}]
^{(p-2)/2p}|z_{j}|_{p}=\omega ([k_{j}+1,k_{j+1}])^{p-2)/2p}|z_{j}|_{p}
$$
and define
$$
z^{*}_{j}(x)=|z_{j}|^{-2}_{2}<z_{j},x>=|z_{j}|_{2}^{-2}\sum^{k_{j+1}}_
{n=k_{j}+1}e^{*}_{n}(z_{j})e^{*}_{n}(x)w^{2}_{n}.$$
applying H\"{o}lder's inequality as above shows that
$$
|z^{*}_{j}(x)z_{j}|_{2}\leq |x_{\left| [k_{j}+1,k_{j+1}]\right. }|_{2}
\text{ and }|z^{*}_{j}(x)z_{j}|_{p}\leq c|x_{\left|
[k_{j}+1,k_{j+1}\right. ]}|_{p}.$$
Combining this observation with the isomorphic classification of subspaces
spanned by block basic sequences, we arrive at a prototype for a complemented
subspace of $X_{p}$  isomorphic to $X_{p,w'}$.

\proclaim{Proposition 1.2} If $(z_{j} )$ is a normalized block basis of 
the natural basis of
$X_{p} ,\; (E_{j} )$ is a sequence of disjoint subsets 
of $\Bbb N    $ and there are positive
constants $c$ and $\delta $ such that for all $j$
\roster
\item"{\rm a)}" $ |z_{j |E_{j}}  |_{2}\geq \delta |z_{j} |_{2}$
\item"{\rm b)}" $ c|z_{j |E_{j}}  |_{2}\geq \omega (E_{j})^{(p-2)/2p}$
\endroster
Then $[z_{j} ]$ is a $\max\{ \delta^{-1}  ,c\} $ complemented subspace of
$X_{p}$  isomorphic to $X_{p,w'},$
where $w^{\prime }_{j} = \omega (E_{j})^{(p-2)/2p}.$

\endproclaim
\demo{Proof} A block basis of the natural basis of $X_{p}$  
spans a subspace isomorphic
to $X_{p,w''}$     where $w^{\prime \prime }_{j} = |z_{j} |_{2} /|z_{j} |_{p}
  = r(z_{j} ).$ By H\"{o}lder's inequality and b), 
$$|z_{j} |_{p}\,\omega(E_{j})^{(p-2)/2p}\geq |z_{j |E_{j}}  |_{p}\,\omega 
(E_{j})^{(p-2)/2p}\geq 
|z_{j|E_{j}}|_{2}\geq \delta
|z_{j}|_{2}.$$
Hence
$$
\delta^{-1}\omega(E_{j} )^{(p-2)/2p}\geq 
r(z_{j})\geq c^{-1}\omega (E_{j} )^{(p-2)/2p}$$
and thus $X_{p,w'}$     is isomorphic to
$X_{p,w''}.$  Define a projection onto $[z_{j}:j \in \Bbb N    ]$ by  
$$\align
               Px &= \sum^{\infty }_{j=1}z^{*}_{j}(x)z_{j}  =
\sum^{\infty }_{j=1}|z_{j|E_{j}}|_{2}^{-2}<z_{j|E_{j} } ,x >z_{j} \\
  &=  \sum^{\infty }_{j=1}|z_{j|E_{j}}|_{2}^{-2}\biggl[
\sum_{n\in E_{j}}e_{n}^{*}(z_{j})e_{n}^{*}(x)w_{n}^{2}\biggr]z_{j} 
\endalign$$
Clearly $P$ is the required operator if it is bounded. The computations above
using b) in the form 
$$
                 (c|z_{j |E_{j}}  |_{p} )|z_{j |E_{j}}  |_{2}\geq \omega 
(E_{j} )^{(p-2)/2p}        |z_{j |E_{j}}  |_{p}
$$
show that
$$\align
||Px|| &\leq \max \{ \biggl[\sum^{\infty }_{j=1}|z^{*}_{j}
(x)z_{j} |^{2}_{2}\biggr]^{1/2},\;\biggl [\sum^{\infty  }_{j=1}|z^{*}_{j}
(x)z_{j} |^{p}_{p}\biggr]^{1/p}\} \\
&= \max \left\{ \biggl[\sum^{\infty }_{j=1}|z^{*}_{j}(x)z_{j
|E_{j}}|^{2}_{2}  | z_{j} |^{2}_{2}/|z_{j|E_{j}} |^{2}_{2}\biggr]^{1/2},\; 
\biggl[
\sum^{\infty }_{j=1}|z^{*}_{j}(x)z_{j |E_{j}}  |^{p}_{p}|z_{j}|_{p}^{p}
/ |z_{j|E_{j}} 
|^{p}_{p}\biggr]^{1/p}\right\}\\
&\leq \max \left\{ \biggl[\sum^{\infty}_{j=1} |x_{|E_{j}}  |^{2}_{2} \delta 
^{-2}\biggr]^{1/2},\; \biggl[\sum^{\infty }_{j=1}  |x_{|E_{j}} 
|^{p}_{p} c^{p}\biggr]^{1/p}\right\} \leq \max \{ \delta ^{-1}  ,c\} ||x||.
\qed
\endalign$$
\enddemo

     Next we will prove our characterization of the complemented subspaces of
$X_{p}$  which contain $X_{p}$  complemented and thus are isomorphic to 
$X_{p}$, \cite{JO2}.

\proclaim{Theorem 1.3} Suppose that $X$ is a complemented subspace of 
$X_{p}.$ Then the
following are equivalent.
\roster
\item"{\rm 1)}" $X$ contains a complemented subspace isomorphic to $X_{p}$
\item"{\rm 2)}" There exist positive constants $c$ and $\delta $ 
such that for every
$\epsilon > 0 $
        there is an $\epsilon',\; 0 < \epsilon' < \epsilon ,$ such that for 
every $N \in  \Bbb N    $ there is an 
        $x \in  X,\; ||x|| = 1$ and a finite set $E \subset  \{ N,\;  N+1,\; 
\dots \} $ such that
\itemitem{\rm a)} $||x_{|[1,N]}|| < N^{-1}$
\itemitem{\rm b)}  $|x_{|E} |_{2}\geq \delta |x|_{2}$
\itemitem{\rm c)}  $\epsilon \geq c|x_{|E} |_{2}\geq \omega (E)^{(p-2)/2p}\geq 
\epsilon '.$
\endroster
\endproclaim

\demo{Proof} Suppose that 2) is satisfied. Let $\epsilon _{k}  =
k^{-1}$.  By 
induction we may choose
for each $k$ a sequence $(x_{k,j})$ of norm one elements of $X$ which are a
perturbation of a block basis of the basis of $X_{p}$  satisfying b) and c)
(for $\epsilon^{\prime }_{k}$ and $E_{k,j})$. Clearly we may assume that
$ E_{k,j}\cap E_{k,m} =\phi $ for $j \neq m$.
By a simple diagonalization argument we can find sets ${\Cal F}_{k}\subset
\Bbb N    $ such that $(x_{k,j})^{\infty }_{k=1,\; j\in {\Cal F}_{k}}$ is 
equivalent to a block basis of $X_{p}$, the sets $E_{k,j},\; k\in 
\Bbb N    ,\; j\in {\Cal F}_{ k}$  are disjoint, and $(\epsilon^{\prime }_{k})
^{2p/(p-2)}\text{card}\;{\Cal F}_{k}\geq 1$ for each $k$. It now follows
from Proposition 1.2 and standard perturbation arguments that
$Y = [x_{k,j}:k \in \Bbb N    , j \in {\Cal F}_{k}]$ is 
isomorphic to $X_{p}$  and that $Y$ is complemented in $X_{p} $.

     For the converse we will actually show that if we take as our isomorph
of $X_{p}$  the special representation $X_{p,w'}$,     where $w'$ is actually 
a doubly indexed sequence $w' = (w_{k,j})$,  where $w_{k,j}= w_{k}$  for all 
$j,\; \lim w_{k} = 0,$ and $\sum^{\infty
}_{k=1}w^{2p/(p-2)}_{k}=\infty $, then
the images of a subsequence of the basis satisfy the properties in 2). Thus we
suppose that $Y$ is a complemented subspace of $X$ and that $T$ is an 
isomorphism of
$X_{p,w'}$ onto $Y.$ By passing to a subsequence of the basis of $X_{p,w'}$ 
and using a
standard perturbation argument, we may assume that $Y$ is the span of a block of
the basis of the containing $X_{p}.$  Let $(y_{i} )$ be the normalized basis 
of $Y$ and let
$F_{i}$  be the support of $y_{i}$  relative to the basis of $X_{p}.$  Let 
$y^{*}_{i}$ denote the biorthogonal functional to $y_{i} .$  Because $Y$ is 
complemented in $X_{p},$ we may assume that each $y^{*}_{i}$ is defined on 
$X_{p}$  and $\sup ||y^{*}_{i}|| \leq ||T^{-1}||\;  ||Q||$ where $Q$ is the
projection onto $Y.$  Because $Y$ is reflexive, $y^{*}_{i}(e_{j})
\rightarrow 0$ as $i \rightarrow \infty $ for each $j$
where $e_{j}$  denotes the $j$th basis vector of $X_{p}.$ Thus we may assume 
by passing to a subsequence, a perturbation argument and perhaps enlarging 
the sets $F_{i}$ slightly that $y^{*}_{i}(x) \neq 0$ only if $x_{|F_{i}}\neq
0.$  In other words $y^{*}_{i}(x) = y^{*}_{i}(x_{|F_{i}}  ).$
Also it follows from this that the projection $Q$ onto $Y$ is given by
$Qx =  \sum^{\infty }_{i=1} y^{*}_{i}(x)y_{i} $.

     Fix $i$ and let $E_{i}  = \{ j \in  F_{i} : |y_{i} (j)| \geq \rho w^{2/
(p-2)}
_{j} |y_{i} |^{-2/(p-2)}_{2} \} $ and assume that b) is satisfied for $y_{i}$
  and $E_{i}.$ 
$$\align 
\omega (E_{i}) =\sum_{j\in E_{i}}w_{j}^{2p/(p-2)} &\leq
\rho^{-2}\sum_{j\in E_{i}}|y_{i}(j)|^{2}|y_{i}|^{4/(p-2)}_{2}w^{2}_{j}\\
&=\rho^{-2}|y_{i|E_{i}}|^{2}_{2}|y_{i}|^{4/(p-2)}_{2}\\
&\leq \rho^{-2}\delta^{-4/(p-2)}|y_{i|E_{i}}|^{2p/(p-2)}_{2}.
\endalign$$

Thus if $\rho $ is independent of $i,$ condition b) will imply the middle 
inequality
in condition c) (with $c = \rho ^{-(p-2)/p}\delta ^{-2/p}).$  Also observe 
that because 
$$|y_{i}|_{p}\,\omega (E_{i} )^{(p-2)/2p}\geq |y_{i |E_{i}} 
|_{p}\,\omega (E_{i})^{(p-2)/2p}\geq  |y_{i_{ |E_{i}} } |_{2}\geq \delta |y_{i} 
|_{2}$$
the third
inequality in c) will be satisfied if $|y_{i} |_{2}$  is bounded away from zero.
Hence it is sufficient to show that for some $\rho  > 0$ there is a $\delta , 
\; 0 < \delta < 1,$
such that if ${\Cal E}_{\delta }  = \{ i:|y_{i |E_{i}}  |_{2}\geq \delta |
y_{i }|_{2} \} $ then for every $\epsilon _{1}  > 0$ there is an
$\epsilon _{2}  > 0$ such that $\epsilon _{1}\geq  |y_{i} |_{2}\geq  \epsilon
_{2}$  for infinitely many $i$ in ${\Cal E}_{\delta } .$  Then c) will be
satisfied with $c\epsilon _{1}  = \epsilon $ and $\epsilon ' = \delta \epsilon
_{2} .$

     Note that      $ \sum_{j\notin
E_{i}} |y_{i} (j)|
^{p}\leq \sum_{j\notin E_{i}} |y_{i} (j)|^{2} \rho ^{p-2}w^{2}_{j} 
|y_{i} |^{-2}_{2}\leq \rho^{p-2}$.   Hence
$||y_{i |E_{i}}||\geq [1 - \rho^{p-2}]^{1/p}$ and $|y_{i |F_{i}\backslash 
E_{i}}|_{p}\leq
\rho ^{1-2/p}.$  Thus if $\rho $ is small,
$[y_{i |E^{ c}_{i}} :i \notin {\Cal E}_{\delta }]$ is not better than $\rho
^{-1+2/p}$       equivalent to $[y_{i} :i\notin
{\Cal E}_{\delta }].$

     For each $K$ define ${\Cal M}_{K}  = \{ i: |y^{*}_{i}(y_{i |E_{i}}  )|\leq 
K |y_{i |E_{i}}  |_{2} /|y_{i} |_{2} \} .$
Observe that because $||y^{*}_{i}|| \leq  ||T^{-1}||\;   ||Q||$,
$$\align
|y^{*}_{i}(y_{i |F_{i}\backslash E_{i}})|&\leq  ||T^{-1}||\;  ||Q||
\max \{ |y_{i |F_{i}\backslash E_{i}} |_{2} ,\; |y_{i|F_{i}\backslash
E_{i}}|_{p}\}\\  
&\leq ||T^{-1}||\;  ||Q||\max  \{ |y_{i}|_{2},\rho^{1-2/p}\}.
\endalign$$
Thus if we consider only those $y_{i}$  with $|y_{i} |_{2}\leq \rho ^{1-2/p}$ 
      we have that 
$$
                  |y^{*}_{i}(y_{i |E_{i}})|\geq 1 - ||T^{-1}||\;  ||Q||\rho
^{1-2/p}. 
$$
Under our assumption on the sequence $(y_{i} ),$ the span of such $y_{i}$  
is still
isomorphic to $X_{p,w}.$  From now on we will assume that $\rho^{1-2/p}\leq 
(||T^{-1}||\;   ||Q||2)^{-1}$   and
thus that $|y^{*}_{i}(y_{i|E_{i}})|\geq 1/2$ for all $i .$ In this way we can 
work with ${\Cal M}_{K}$
instead of ${\Cal E}_{\delta }$  since for such $i$ if $i \in {\Cal M}_{K}$
then $i \in {\Cal E}_{1/ 2K}.$

     Let us now see how the projection onto $[y_{i}]$ acts on the span of
$[y_{i |E_{i}}].$
Our assumptions on the $y_{i}$'s imply that $ Q
\sum^{\infty }_{i=1} a_{i} y_{i|E_{i}}   =\sum^{\infty }_{i=1}
 a_{i} y^{*}_{i}(y_{i |E_{i}}  )y_{i} $.  Hence 
$$\align
   \biggl| Q\sum^{\infty }_{i=1} a_{i} y_{i |E_{i}}\biggr| _{2}&=\biggl|
\sum^{\infty }_{i=1}a_{i}y^{*}_{i}(y_{i |E_{i}}  )y_{i}\biggr| _{2}\\
&\geq 
K\biggl[ \sum_{i\notin {\Cal M}_{K}} |a_{i} |^{2}| 
y_{i|E_{i}}| ^{2}_{2}| y_{i} | ^{-2}_{2}| y_{i}|
^{2}_{2}\biggr] ^{1/2}     
                       =    K\biggl[ \sum_{i\notin{\Cal M}_{K}}|
a_{i} |^{2} |y_{i |E_{i}}  |^{2}_{2}\biggr] ^{1/2}      .
\endalign$$
If $a_{i} = |y_{i |E_{i}}  |_{2}^{2/(p-2)},$ for $i \notin
{\Cal M}_{K}$ 
and $i \leq N,$ and 0 else, then

$$\align
\biggl |\biggl|  \sum^{N}_{i=1} a_{i} y_{i |E_{i}}\biggr|\biggr| 
&= \max \biggl\{ \biggl[ 
\sum_{i\notin {\Cal M}_{K}}|a_{i}|^{p}|y_{i |E_{i}} 
|_{p}^{p}\biggr] ^{1/p},\; \biggl[ \sum_{i\notin {\Cal M}_{K}}
|a_{i} |^{2}
 |y_{i |E_{i}}  |^{2}_{2}\biggr] ^{1/2}\biggr\}\\ 
&\leq \max\biggl\{ \biggl[ \sum_{i\notin {\Cal M}_{K}}|y_{i|E_{i}}
|^{2p/(p-2)}_{2}|y_{i}|^{p}_{p}\biggr]^{1/p},\; \biggl[ \sum_{i\notin
{\Cal M}_{K}}|y_{i|E_{i}}|_{2}^{2p/(p-2)}\biggr]^{1/2}\biggr\} \\
&\leq \max \biggl\{ \biggl[ \sum_{i\notin {\Cal
M}_{K}}|y_{i|E_{i}}|^{2p/(p-2)}_{2}\biggr] ^{1/p},\; \biggl[ 
\sum_{i\notin {\Cal M}_{K}}|y_{i|E_{i}}|_{2}^{2p/(p-2)}\biggr] 
^{1/2}\biggr\} \\
&= \sum_{i\notin {\Cal M}_{K}}\biggl[
|y_{i|E_{i}}|_{2}^{2p/(p-2)}\biggr] 
^{1/2}, \qquad \text{if }\sum_{i\notin {\Cal
M}_{K}}|y_{i|E_{i}}|^{2p/(p-2)}_{2}\geq 1.
\endalign$$
This implies that 
$$
||Q||\biggl[ \sum_{i\notin {\Cal
M}_{K}}|y_{i|E_{i}}|_{2}^{2p/(p-2)}\biggr]^{1/2}\geq K\biggl[ 
\sum_{i\notin {\Cal M}_{K}}|y_{i|E_{i}}|_{2}^{2p/(p-2)}\biggr] 
^{1/2},
$$
if $ \sum_{i\notin {\Cal
M}_{K}}|y_{i|E_{i}}|_{2}^{2p/(p-2)}\geq 1.$   Therefore, if $K \geq | |Q||,$
$$ \sum_{i\notin {\Cal M}_{K}}|y_{i|E_{i}}|_{2}^{2p/(p-2)}\leq 1.$$
Because $ |y_{i |E_{i}}  |_{p}\geq [1 - \rho ^{p-2}]
^{1/p},$ ( We are assuming that $\rho < 1.$)
this implies that $ [y_{i|E_{i}}||y_{i |E_{i}}||^{-1}  :i 
\notin {\Cal
M}_{K}]$ is equivalent to the basis
of $\ell_{p} .$  Therefore for any $\epsilon '$ small enough only finitely 
many of the $y_{i}$'s
with $|y_{i} |_{2}\geq \epsilon '$ have index not in ${\Cal M}_{K} .$  
Indeed, if this were not the case,
then there would be a subsequence of $(y_{i} )$, say $(y_{i} )_{i\in
M}$,    
such that $M \subset {\Cal M}^{c}_{K}$  and
$\epsilon  > 0 , |y_{i} |_{2}\geq \epsilon $  for all $i \in M.$  However 
this would imply that $Q$ is an
isomorphism from $[y_{i|E_{i}}:i \in M]$, which is isomorphic to 
$\ell _{p},$ onto $[y_{i} :i \in M],$
which is isomorphic to $\ell _{2}.$

     It follows that $(y_{i} )_{i\in {\Cal M}_{K}}$     is equivalent to the 
basis of $X_{p,w'},$ where $w' = (w^{\prime }_{i})$ and for each $\epsilon > 0
$
there is a $\epsilon ' > 0$ such that $\epsilon > w^{\prime }_{i}\geq \epsilon
',$ for infinitely many $i.$  Note that because $|y^{*}_{i}(y_{i|E_{i}})|\geq 
1/2 ,$ for any $i \in {\Cal M}_{K},$\ $|y_{i|E_{i}} |_{2}\geq \delta  |y_{i}
|_{2}  ,$ where $\delta = 1/2K ,$ i.e., $i \in {\Cal E}_{\delta }.$ 
\qed
\enddemo

{\smc Remark 1.4:} If $X$ is not isomorphic to $X_{p}$  
we  can use part of the proof that 1)
implies 2) to get a natural way of splitting vectors in $X$ into a piece with
large ratio and a piece with small ratio. Indeed suppose $X$ is a complemented
subspace of $X_{p}$  and that $P$ is the projection onto $X.$  By 
\cite{JJ} or \cite{JO2} we may
assume that $P$ is bounded in the norm $|\cdot |_{2}$  as well. Suppose that 
2) fails for
$\delta < |P|^{-1}_{2} ,\; c$ and $\epsilon .$  Choose positive constants 
$\epsilon
',\; \rho ,\; \alpha ,$ and $\beta $ such that 
$$\align
 \epsilon ' &<\min \{ \epsilon ,\delta \alpha \} ,\\
\beta &< \min \{ (1 - \delta |P|_{2} )/(||P||) ,\; \epsilon /c\} ,\\
 \rho &\leq \min \{ c^{-p/(p-2)}\delta ^{2/(p-2)},\; \beta^{p/(p-2)}\} ,\\
 \beta >\alpha &\geq \max \{ \beta \delta |P|_{2}/(1-\beta ||P||),\; \beta
^{2}||P||/(1-\delta |P|_{2}).
\endalign$$
Let $N$ be an integer so that a), b), and c) of 2) fail for $\epsilon '$ and
$N$ and suppose that $x \in X ,\; x|_{[1,N]}=0,\alpha < r (x) < \beta $ 
and $||x|| = 1.$  Let
$E_{x}   = \{ j : |x(j)| \geq \rho w^{2/(p-2)}_{j}|x|^{-2/(p-2)}_{2}\} .$  
As in the proof above the 
choice of
$\rho $ guarantees that the middle inequality in 2) c) is satisfied by
$x_{|E_{x}}.$  Because $r(x) < \beta \leq \epsilon /c,$ the first inequality 
in 2) c) is also satisfied.
Finally if $|x_{|E_{x}}  |_{2}\geq \delta |x|_{2} ,\; \omega (E_{x}
)^{(p-2)/2p}\geq |x_{|E_{x}}|_{2}\geq \delta |x|_{2}\geq \delta \alpha 
\geq \epsilon '$ and thus
all of the inequalities in c) are satisfied. The failure of 2) then implies
that $|x_{|E_{x}}  |_{2}  < \delta |x|_{2} .$

     Let $y = P(x_{|E_{x}}  )$ and $ z = x - y = 
P(x_{|E^{c}_{x}}  ).$  We 
claim that $r(y)\leq \alpha $ and
$r(z)\geq \beta .$  Indeed, 
$$\align
|y|_{2}&\leq |P|_{2} |x_{|E_{x}}  |_{2}  < |P|_{2}
\delta |x|_{2}  \leq |P|_{2}\delta \beta
\\ &\leq \alpha (1-\beta ||P||)\leq \alpha
(|x|_{p}-||P||\; ||x_{|E^{c}_{x}}||)\leq \alpha |y|_{p}
\endalign$$
since $|x_{|E_{x}^{c}}|
_{p}\leq \rho ^{(p-2)/p}\leq \beta $ and $|x|_{2}\leq \beta $, and
$$|z|_{2}\geq (1-\delta |P|_{2})\alpha \geq (1-\delta |P|_{2})\beta
^{2}||P||/(1-\delta |P|_{2})\geq \beta ||P||\; ||x_{|E_{x}^{c}}||\geq \beta
|z|_{p}.$$

     Thus any $x \in  X$ with support in $\{ N+1,N+2,\dots \} $ can be split 
into an
element with ratio greater than $\beta $ and one with ratio smaller than
$\alpha .$  If this
could be accomplished in a linear fashion it would follow that $X$ is then
isomorphic to a complemented subspace of $\ell  _{p}\oplus \ell _{2}.$

\newpage
\heading
2. Complemented subspaces of $X_{p}$  which are isomorphic to $\ell 
_{p}\oplus \ell _{2} $ \endheading

     In this section we look at some ways of discriminating between
complemented subspaces of $X_{p}$  
which are isomorphic to complemented subspaces
of $\ell_{p}\oplus \ell_{2}$  and those isomorphic to $X_{p} .$  
First we will examine how the
conditions in Theorem 1.3 fail if $X$ is isomorphic to $\ell_{p}\oplus
\ell_{2} .$  Below $P_{n}$ denotes the basis projection onto the span of the
first $n$ elements of the basis of $X_{p}.$

\proclaim{Proposition 2.1}  Suppose that $Z,\; X,\; U,$ and $W$ are subspaces of
$X_{p}$  such that
$Z \subset X = U \oplus W,\; U$ is isomorphic to $\ell_{2}$  and $W$ 
is isomorphic to $\ell_{p}.$ Suppose that
$Z$ has a normalized $K$ unconditional basis $(z_{n} ).$  Let $
\beta  = 
\lim_{n\to \infty } r(z_{n} )$ and $\beta ' =
\lim_{n\to \infty } \inf \{ b:\text{for every }\epsilon >0\text{ there exists 
}
u \in U \text{ such that }||P_{n} u||<\epsilon \text{ and }r(u) \leq b\} .$ If
$\beta > 0,\; \beta ' \leq 1,$ and $P$ is a projection from $X$ onto $Z$ 
then $||P|| \geq \beta '/K\beta .$
\endproclaim
\demo{Proof} Let $z_{n}  = u_{n}  + w_{n}$  where $u_{n}\in U$ 
and $w_{n}  \in  W.$ 
By passing to subsequences
and a standard perturbation argument we may assume that $(u_{n} )$ and  
$(w_{n} )$ are
block bases of the basis of $X_{p} .$ (It could happen that
$||w_{n}||\to 0,$ but then $\beta ' \leq \beta .)$  Moreover we may 
assume that the projection $P$ composed with the
corresponding basis projection $Q$ acts disjointly with respect to the
subsequence $(z_{n} )_{n\in M},$ i.e., $QP$ is a projection onto $[z_{n}:
n\in M]$ and $QPu_{n}  = \tau_{n} z_{n}$
and $QPw_{n}  = (1 - \tau_{n} )z_{n}.$  Because $(w_{n})$ is equivalent to 
the usual unit vector
basis of $\ell_{p},\; p > 2,\; \beta > 0,$ and $(z_{n} )$ is equivalent to 
the usual unit vector
basis of $\ell_{2},$ it follows that $\tau_{n}\to 1.$

Because $W$ is isomorphic to $\ell_{p},\;|w_{n} |_{2}
\to 0$ and thus $|u_{n} |_{2} - |z_{n} |_{2} \to 0.$
Therefore $$\align
\limsup ||u_{n}|| = \limsup \max\{ |u_{n} |_{2} ,\; |u_{n}
|_{p}\} &\leq \lim \sup \max \{ \beta ,\; |u_{n} |_{2} /r(u_{n} )\}\\ 
&\leq \max \{ \beta ,\beta /\beta '\} =\beta /\beta '.
\endalign$$ 
Consequently 
$K||P||\beta /\beta '\geq \limsup ||Q||\; ||P||\; ||u_{n}||\geq 1.$
\qed
\enddemo
\proclaim{Corollary 2.2}Suppose that $X,\; U,$ and $W$ are subspaces of 
$X_{p}$  
which satisfy the hypotheses of Proposition 2.1 and $X$ is complemented in 
$X_{p}$  with projection $P.$
Then for any $c$ and $\delta $ and $\epsilon < \beta 'c\delta /\max\{
c,\delta^{-1}\} ,$  there is no $\epsilon ',\; 0 < \epsilon '  < \epsilon ,$
such that for every $N \in \Bbb N    $ there is an $x \in  X,\; ||x|| = 1$ 
and a finite set $E \subset \{ N,\; N+1,\; \dots \} $ such that
\roster
\item"{\rm a)}" $||x_{|[1,N]}||< N^{-1}$
\item"{\rm b)}" $|x_{|E} |_{2} \geq \delta |x|_{2}$
\item"{\rm c)}" $\epsilon \geq c|x_{|E} |_{2} \geq 
\omega (E)^{(p-2)/2p}\geq \epsilon '.$
\endroster
\endproclaim
\demo{Proof} Suppose $\epsilon '$ exists for some $c,\; 
\delta ,$ and $\epsilon .$ Then there is a sequence
$(z_{n} )$ of norm one vectors in $X$ which is a perturbation of 
a block basis of the
basis of $X_{p}$  and disjoint sets $(E_{n} )$ such that for all $n$
\roster
\item"{\rm b)}" $|z_{n|E_{n}} |_{2} \geq \delta |z_{n}|_{2}$;
\item"{\rm c)}" $\epsilon \geq c|z_{n|E_{n}} |_{2} 
\geq \omega (E_{n})^{(p-2)/2p}\geq \epsilon
'.$
\endroster
Then $[z_{n} :n \in \Bbb N    ]$ is complemented in $X_{p}$  by a projection 
of norm at most
$\max \{ c,\delta^{-1}\}$  and $ r(z_{n} ) \leq \epsilon /c\delta .$ 
Thus by 
the previous proposition $||P|| \geq \beta 'c\delta /\epsilon $
and hence  $$\epsilon \geq \beta 'c\delta /\max \{ c,\delta^{-1}\} .
\qed$$
\enddemo

     We now turn our attention to the classification of the complemented
subspaces of $X_{p}.$ It was shown in \cite{JO2} that if a 
complemented subspace of 
$X_{p}$ has an unconditional basis then it is isomorphic to $\ell_{p},\;
\ell_{2},\; \ell_{p}\oplus \ell_{2},$ or $X_{p}.$
In \cite{AC} the same conclusion was established if $X$ has 
a ``$p,2$ F.D.D." 
Thus it seems likely that same result holds without the additional 
assumptions. We will next look at some well known results but recast in terms 
of the ratio of the 2-norm and $p$-norm.

     To begin let us recall that results of Kadec and Pelczynski
\cite{KP} give a natural criterion for isomorphs of $\ell_{2}$  contained 
in $X_{p},
p > 2,$ namely, a subspace $X$ of $X_{p}$  is isomorphic to $\ell_{2}$  
if and only if there is a constant $C > 0$ such that $r(x) = |x|_{2}
/|x|_{p}\geq C$ for all $x \in X,$ i.e., $h(X) \geq C.$
It follows from \cite{JO1} that if a complemented subspace of $X_{p}$  does not 
contain $\ell_{2}$  then it is isomorphic to $\ell_{p} .$  A standard gliding 
hump argument yields the
following criterion. (Below $Q_{N}$  denotes the projection onto the span of the
basis vectors of $X_{p}$  with index greater than $N.$)

\proclaim{Proposition 2.3}  A complemented subspace $X$ of 
$X_{p}$  is isomorphic to
$\ell_{p}$  if and
only if for every $\epsilon > 0$ there is an $N \in \Bbb N    $ such that if 
$x \in Q_{N} X$ then $r(x) < \epsilon .$
\endproclaim

     Theorem 1.3 gives a criterion for identifying complemented subspaces
isomorphic to $X_{p}$, however it seems to be rather difficult to formulate
useful conditions which identify complemented subspaces isomorphic to $\ell_{2}
\oplus \ell_{p}$ or a complemented subspace of it. Here are some attempts 
at such criteria.
\proclaim{Proposition 2.4}  Let $X$ be a complemented subspace of $X_{p}$  and 
suppose that $Z$ is
a subspace of $X$ and $\epsilon ,\;\beta ,$ and $\beta '$ are positive 
constants with $\epsilon \leq 1$ such that
\roster
\item"{\rm a)}" for all $z \in Z,\; r(z) \geq \beta '$
\item"{\rm b)}" if $x \in  X$ and $r(x) > \beta $ 
then there exists $z \in Z$ such 
that $|x - z|_{2}  < \epsilon |x|_{2} .$
\endroster
Then $X$ is isomorphic to a complemented subspace of $\ell_{p}\oplus
\ell_{2}$ and conversely.
\endproclaim

\demo{Proof} Condition a) implies that $Z$ is isomorphic to $\ell_{2} .$  
Let $Y$ 
be the kernel
of the orthogonal projection from $X$ onto $Z.$  If $y \in Y$ and $r(y) > \beta
$ then by b) there exists a $z \in Z$ such that $|y - z|_{2}  < \epsilon
|y|_{2}.$ But $y$ is orthogonal to $z$ so
we have that $|y|_{2}\leq [|y|^{2}_{2}  + |z|^{2}_{2}]^{1/2}=|y-z|_{2}  < 
\epsilon 
|y|_{2}$, an impossibility.
Therefore $r(y) \leq \beta $ for all $y \in Y$ and $Y$ is isomorphic to
$\ell_{p} .$

     If $X$ is isomorphic to $\ell_{2} $ or $\ell_{p}$  the converse follows 
easily from our
earlier observations.  Thus by the results of Edelstein and Wojtasczyk 
\cite{EW} we
may assume that $X$ is isomorphic to $\ell_{p}\oplus \ell_{2}$  and let $U$
 and $W$ be the
complementary subspaces with $U$ isomorphic to $\ell_{2}$  and $W$ isomorphic 
to $\ell_{p}$.
Because $U$ is isomorphic to $\ell_{2}$  there is a constant $1 \geq  \beta ' 
> 0$ such that for all
$u \in U,\; r(u) \geq \beta '.$  We may also assume that $W$ is the 
kernel of the orthogonal
projection $Q$ onto $U.$

     Then for any $x \in X,$ 
$$\align
   |x - Qx|_{2}&\leq  r(W)||x - Qx|| \leq r(W)(1 + ||Q||)||x||\\
&\leq r(W)(1 + 1/\beta ')|x|_{2} \max\{ 1,1/r(x)\} .
\endalign$$
Because $W$ is isomorphic to $\ell_{p}$  it has a basis, let $R_{n}$  denote 
$I - Q$ composed
with the basis projection onto the span of the first $n$ elements of the basis
of $W.$ Let $K = \sup ||R_{n} ||.$  Choose $n$ so large that if $Y = (I - R_{n}
 )W,\; r(Y) < (2(1 + 1/\beta ')(1 + K))^{-1}$. Then if $x \in Y + Z$ the 
above computation shows
that b) is satisfied with $\epsilon = 1/2.$  Because $R_{n}$  is finite rank 
there exists a
$\beta > 1$ such that if $r(x) > \beta $ then $||R_{n} x|| < ||x||/4(1+||Q||).$
 Now if $r(x) > \beta ,$ 
$$\align
    |x - Qx|_{2}  &\leq (1 + ||Q||)||R_{n} x|| + |(I - R_{n} )x - Q(I - R_{n}
 )x|_{2}\\
                  &\leq (1 + ||Q||)||x||/4(1 + 1/\beta ') + r(Y)||(I - R_{n} )
x - Q(I - R_{n} )x||\\
                  &\leq |x|_{2}\max \{ 1,1/r(x)\} /4 + r(Y)(1 + ||Q||)||I -
R_{n}||\; ||x||\\
                  &\leq |x|_{2}\max \{ 1,1/r(x)\} /4 + r(Y)(1 + 1/\beta ')
(1 + K)|x|_{2} \max\{ 1,1/r(x)\} .\\
                  &\leq (3/4)|x|_{2} .
\qed
\endalign$$
\enddemo

{\smc Remark 2.5:}  Condition b) may be replaced by
\roster
\item"{\rm b$'$)}" if $x \in  X$ and $r(x) > \beta ''$ then there 
exists $z \in Z$ such that $|x - z|_{2}  < \epsilon ||x||.$
\endroster
To see this note that if b$'$) holds then b) holds with $\beta  = \max \{ 
\beta '',\epsilon \} $ and $\epsilon = 1.$

\bigpagebreak
    To get a similar theorem but with the hypothesis on the $\ell_{p}$  part 
we seem
to need to assume the existence of a projection.
\proclaim{Proposition 2.6}  Let $X$ be a complemented subspace of $X_{p}$  and 
suppose that $Y$ is
the range of a projection $P$ on $X$ and $\epsilon ,\; \alpha ,$ and $\alpha '$
 are positive constants with
$\epsilon \leq  ||I - P||^{-1}$   such that
\roster
\item"{\rm a)}" for all $y \in Y,\; r(y) \leq \alpha '$
\item"{\rm b)}" if $x \in X$ and $r(x) < \alpha $ then there 
exists $y \in Y$ such that $||x - y|| < \epsilon ||x||.$
\endroster
Then $X$ is isomorphic to a complemented subspace of $\ell_{p}\oplus
\ell_{2}$  and conversely.
\endproclaim
\demo{Proof} Condition a) implies that $Y$ is isomorphic to $\ell_{p} .$  
Let $Z$ 
be the kernel
of the projection $P$ from $X$ onto $Y.$  If $z \in Z$ and $r(z) < \alpha $ 
then by b) there exists a $y \in Y$ such that $||z - y|| < \epsilon ||z||.$ 
But $Pz = 0$ and $Py = y$ so
we have that $||z|| = ||(I - P)(z - y)|| < ||I - P||\epsilon ||z|| \leq ||z||,$
 an impossibility.
Therefore $r(z) \geq \alpha $ for all $z \in Z$ and $Z$ is isomorphic to
$\ell_{2} .$
   
  As in the proof of Proposition 2.6 the converse easily reduces to the
case that $X$ is isomorphic to $\ell_{2}\oplus \ell_{p}.$  So we again let 
$U$ and $W$ be the
complementary subspaces with $U$ isomorphic to $\ell_{2}$  and $W$ isomorphic 
to $\ell_{p}$ and let
$\beta '$ be a constant such that $1 \geq \beta ' > 0$ and for all $u \in U, 
r(u) \geq \beta '.$ As before
we will assume that $W$ is the kernel of the orthogonal projection $Q$ onto
$U.$

     Then for any $x \in  X,$ 
$$\align
     ||Qx|| = \max \{ |Qx|_{2} ,|Qx|_{p} \} &\leq \max \{ |Qx|_{2} ,|Qx|_{2} 
/\beta '\}  \leq <Qx,x>^{1/2}   /\beta '\\
&\leq [|Qx|_{2} r(x)||x||]^{1/2}   /\beta ' \leq ||Q||^{1/2}   r(x)^{1/2}
   ||x||/\beta ' .
\endalign$$
Thus if $r(x) < \alpha  = {\beta '}^{2} /||Q||^{3},$\  
$||x - (I - Q)x|| = ||Qx|| < 
||x||/||I - (I - Q)||.$  Because $Y$ is isomorphic to $\ell_{p}$ there is some
$ 
\alpha '$ such that $r(y)\leq \alpha '$ for all $y\in Y.$
\qed
\enddemo

{\smc Remark 2.7:} Propositions 2.3, 2.4, and 2.6 do not really use the 
structure of
$X_{p}$  and thus can be restated for complemented subspaces of $L_{p} $.

\bigpagebreak
     Proposition 2.6 should be compared to the following result for $X_{p}$  
itself.
\proclaim{Proposition 2.8} There does not exist a subspace $Y$ of $X_{p}$  
and positive
constants $\epsilon $ and $\alpha ,\; \epsilon < 1$ such that
\roster
\item"{\rm a)}" $r(Y) < \infty $
\item"{\rm b)}" if $x \in X_{p}$  and $r(x) < \alpha $ then there exists a 
$y \in Y $
with $||x - y||<\epsilon ||x||.$
\endroster
\endproclaim
\demo{Proof} Suppose such a subspace exists. Then there is a normalized 
block basic
sequence $(x_{n} )$ of the $X_{p}$  basis such that $\alpha > r(x_{n} ) >
\alpha /2$ for all $n$ and such
that $[x : n \in \Bbb N    ]$ is norm one complemented in $X_{p}$ with
projection $P.$  By b) 
for each $n$ there is an
element $y_{n}$  of $Y$ such that $||x_{n}  - y_{n}|| < \epsilon
||x_{n}||.$ Because $P$ is norm 1,
$||Py_{n}  - x_{n}|| < \epsilon < 1.$ Hence $||Py_{n}|| >
1 - \epsilon .$  By passing to a subsequence we may assume that $(y_{n} )$ 
is equivalent to
the usual unit vector basis of $\ell_{p}$  and that $(Py_{n} )$ is equivalent 
to a block
basic sequence in $[x_{n} :n \in \Bbb N    ].$  But $(x_{n} )$ is 
equivalent to the unit vector basis
of $\ell_{2}$  and hence so is $(Py_{n} ).$ Because $p > 2$ this is a 
contradiction.
\qed
\enddemo

{\smc Remark 2.9:} The above proposition fails if $\epsilon = 1.$ In this
case the span of a perturbation of a natural basic sequence equivalent
to the basis for $\ell_{p}$ may be used for $Y.$

\enddocument